\documentclass[11pt,a4paper]{article}
\usepackage{epsf,algorithm,epsfig,amsfonts,amsgen,amsmath,amstext,amsbsy,amsopn,amsthm,lineno}
\usepackage{color}
\usepackage{tikz}
\usepackage{cite}
\usepackage{enumerate}
\usepackage{float}
\usepackage{epstopdf}

\usepackage{subcaption}
\usepackage{hyperref}
\usepackage{fancyhdr}

\usepackage{multicol}
\usepackage{amssymb}
\usetikzlibrary{decorations.pathreplacing,calligraphy}
\usetikzlibrary{calc}

\newtheorem{theorem}{Theorem}

\newtheorem{problem}{Problem}

\newtheorem{case}{Case}

\newtheorem{claim}{Claim}

\newcounter{mathitem}

\usetikzlibrary{decorations.markings}
\tikzstyle{vertex}=[circle, draw, inner sep=0pt, minimum size=5pt]

\begin{document}

\title{\Large  On saturation problems for matchings with regularity constraints} 
\date{}
\author{Gang Yang$^{a}$, Zixuan Yang$^{b,c}$\thanks{Corresponding author.}, Shenggui Zhang$^{b,c}$ ~\\[2mm]
\small $^{a}$Graduate School of Environment and Information Sciences, \\
\small Yokohama National University, 79-2 Tokiwadai, Hodogaya-ku, Yokohama 240-8501, Japan\\
\small $^{b}$School of Mathematics and Statistics, \\
\small Northwestern Polytechnical University, Xi'an, Shaanxi 710129, P.R. China\\
\small $^{c}$Xi'an-Budapest Joint Research Center for Combinatorics, \\
\small Northwestern Polytechnical University, Xi'an, Shaanxi 710129, P.R. China\\
}
\maketitle

\begin{abstract}
A graph $G$ is  $F$-saturated if $G$ is $F$-free but for any edge $e$ in the complement of $G$ the graph $G + e$ contains  $F$.  Gerbner et al. (Discrete Math., 345 (2022), 112921)   initiated the study of  $rsat(n,F)$, the minimum number of edges in a regular $n$-vertex $F$-saturated graph, and they posed the problem of for which graphs  $rsat(n, F )$ exists.  Regarding this problem,  we obtain the precise value of $rsat(n,(m+1)K_2)$ for all possible cases, where $(m+1)K_2$ denotes a matching of size $m+1$.
As a natural counterpart, we also determine  the maximum number of edges in a regular $n$-vertex  $(m+1)K_2$-free graph for all $m\ge 1$ and $n\ge 2m+2$.

\smallskip
\end{abstract}

\setcounter{footnote}{0}
\renewcommand{\thefootnote}{}
\footnotetext{E-mail addresses:  {\tt gangyang98@outlook.com (G. Yang), yangzixuan@nwpu.edu.cn (Z. Yang), sgzhang@nwpu.edu.cn (S. Zhang)}}


\section {\large Introduction}

\noindent  Throughout this paper all graphs are finite, undirected, and simple.   For terminology and notations not defined here, we refer the reader to Bondy and Murty\cite{Bon}. 

Let $G$ be a graph with vertex set $V(G)$ and edge set $E(G)$.  The number of vertices of $G$ is called its \emph{order}; and the number of edges in $G$, denoted by $e(G)$, is called its \emph{size}. For a vertex $u\in V(G)$, the \emph{degree} of $u$ in $G$ is denoted by $d_G(u)$. For two disjoint subsets $S,T\subseteq V(G)$, let $E_{G}(S,T)$ denote the set of edges of $G$ joining $S$ to $T$, and $e_G(S,T):=|E_G(S,T)|$.  For $S\subseteq V(G)$, we use $G-S$  to denote the subgraph $G[V (G)\setminus S]$. For two graphs $G_1$ and $G_2$, we use $G_1\cup G_2$ to denote the graph with vertex set $V(G_1)\cup V(G_2)$ and edge set $E(G_1)\cup E(G_2)$.

A \emph{matching} of a graph $G$ is a  subset of $E(G)$ such that no two edges share a vertex in common.
A matching $M$ is a \emph{maximum matching} of a graph $G$ if there does not exist a matching $M'$ in $G$ such that $|M'|>|M|$. The cardinality of a maximum matching of a graph $G$,  denoted by $\nu(G)$, is called the \emph{matching number} of $G$. A \emph{perfect matching} of a graph is a matching covering all vertices. 
A graph $G$  is said to be \emph{factor-critical} if $G-\{v\}$  has a perfect matching for every $v\in V(G)$.

Given two positive integers $p$ and $q$, let $K_{p,q}$ denote the complete bipartite graph with $p$ and $q$ vertices in the two bipartite sets respectively; we use $K_n$ and $C_n$ to denote the complete graph and the cycle of $n$ vertices, respectively.   A complete subgraph is called a \emph{clique}. We write $k$ disjoint copies of $K_n$ as $kK_n$.

For a graph $F$, the widely investigated \emph{Tur\'{a}n number} $ex(n,F)$ are concerned with the maximum number of edges that an $n$-vertex graph can have, if it does not contain $F$ as a subgraph (namely  \emph{$F$-free}).  A natural counterpart is the \emph{saturation number} $sat(n,F)$, the minimum number of edges in $n$-vertex $F$-saturated graphs, where a graph $G$ is called \emph{$F$-saturated} if it is $F$-free, but adding any edge to $G$ creates a copy of $F$.  
In efforts to further understand saturation numbers and the structure of saturated graphs, there has been a substantial amount of research on adding degree conditions to the saturated graphs. For instance, imposing a maximum degree or minimum degree condition has been investigated by Alon, Erd\H{o}s, Holzman, and Krivelevich \cite{Alo1}, Day \cite{Day}, Duffus    and  Hanson \cite{Duf}, F\"uredi and Seress \cite{Fur}, Hanson and Seyffarth \cite{Han}.  In this paper, we study a related function, where one restricts to regular graphs. 

Gerbner,  Patk\'{o}s, Tuza, and Vizer \cite{Ger3}  initiated the study of  the existence of regular $F$-saturated graphs and their number of edges.   Let $rsat(n, F )$ denote the minimum number of edges in a regular $n$-vertex $F$-saturated graph (if such a graph exists).  They proved that $rsat(n, K_3)$ exists for all sufficiently large $n$, and they also  posed the problem (see Problems 4.4 and 4.5 in \cite{Ger3}) in the following.

\begin{problem}[Gerbner,  Patk\'{o}s, Tuza, and Vizer, 2022]
For which graphs does $rsat(n, F )$ exist for some $n$ or  infinitely many values of $n$?
\end{problem}

Regarding this problem,  Timmons \cite{Tim} proved $rsat(n,K_4)$ and $rsat(n,K_5)$ exist for   $n>59$. Davini and Timmons \cite{Dav} proved that for any even integer $l \ge 4$, there is an $n$-vertex regular $C_{l+1}$-saturated graph for all $n > 12l^2 +36l+24$. The  first aim here is to prove $rsat(n,(m+1)K_2)$ exists for some $n$ and determine its value.

 \begin{theorem}\label{main-result1}
Let $n$ and $m\ge 1$ be  two integers.  If $2m+2\le n\le 3m$ and $(n - 2m)\mid m$, 
  \begin{align*}
rsat(n,(m+1)K_2)= \frac{nm}{n-2m};
\end{align*}
otherwise, $rsat(n,(m+1)K_2)$ does not exist.
\end{theorem}

 Gerbner,  Patk\'{o}s, Tuza, and Vizer \cite{Ger1} introduced the  \emph{regular Tur\'{a}n number} $rex(n, F )$,  the maximum number of edges in a regular $n$-vertex $F$-free graph. 
Most of the previous work on regular Tur\'{a}n numbers is given by the extensive study of cages (see \cite{Exo} for a survey), where one forbids all cycles up to a fixed length.  For other graphs, Gerbner,  Patk\'{o}s, Tuza, and Vizer  \cite{Ger1} showed that for non-bipartite $F$ with odd girth $g$, one has $rex(n, F ) \geq n^2/(g + 6)-O(n)$ and asked to determine $\lim \inf_{n\to \infty} rex(n, F )/n^2$ for non-bipartite $F$. This problem was solved independently in \cite{Cam} and \cite{Car2} for graphs $F$ with chromatic number at least 4, and the authors obtained asymptotic results for several graphs with chromatic number 3. Exact results on regular Tur\'{a}n numbers of trees were determined in \cite{Ger2}, and complete graphs were obtained in \cite{Car2, Ger2}. Tait and Timmons \cite{Tait} gave lower bounds on $rex(n, F )$, that are best possible up to a constant factor, when $F$ is one of complete bipartite graphs $K_{2,p}, K_{3,3}$, or $K_{m,p}$ for $p > m!$.
The  second aim  in this paper is   to  determine the regular Tur\'{a}n number of matchings.

 \begin{theorem}\label{main-result2}
Let $n$ and $m$ be   integers such that $n\geq 2m+2$ and $m\geq 1$. Then
  \begin{align*}
rex(n,(m+1)K_2)= \begin{cases}
n(\lfloor\frac{n}{n-2m}\rfloor-1)/2,& \text{if  $\lfloor\frac{n}{n-2m}\rfloor$ is odd;}\\
n(\lfloor\frac{n}{n-2m}\rfloor-2)/2,& \text{if $\lfloor\frac{n}{n-2m}\rfloor$ is even.}
\end{cases}
\end{align*}
\end{theorem}


\noindent\textbf{Remark 1.} There exist several extremal graphs satisfying conditions to attain the
maximum number of edges in Theorem \ref{main-result2}. 
Let  $2m+2\le n$,  and  let  $r\in \{\lfloor n/n-2m\rfloor-1, \lfloor n/n-2m\rfloor-2\}$ be an even integer and $r\ge 2$.  Note that 
\begin{align*}
(n-2m)(r+1)\le (n-2m)\Big(\frac{n}{n-2m}\Big)= n.
\end{align*}
Let $K_{2r_1+1},K_{2r_2+1},\ldots, K_{2r_{n-2m}+1}$ be some  complete graphs such that 
\begin{align*}
\min\{2r_1+1,2r_2+1,\ldots,2r_{n-2m}+1\}\ge r+1,
\end{align*}
 and $n=\sum_{i=1}^{n-2m} 2r_i+1$.   
 For each $i\in \{1,2,\ldots,n-2m\}$, $K_{2r_i+1}$ can be decomposed into $r_i$ edge-disjoint Hamilton cycles, saying $C_{2r_i+1}^1, C_{2r_i+1}^2,\ldots, C_{2r_i+1}^{r_i}$, let $G_i=\bigcup_{j=1}^{r/2}C_{2r_i+1}^{j}$. Take $G=G_1\cup G_2\cup \cdots \cup G_{n-2m}$ as the extremal graph.   One can check that  $G$ is an $r$-regular graph of order $n$ with 
 \begin{align*}
 \nu(G)=\sum_{i=1}^{n-2m}\nu(G_i)=\sum_{i=1}^{n-2m}r_i=m~~(since~n=\sum_{i=1}^{n-2m} 2r_i+1),
  \end{align*}
and
\begin{align*}
e(G)=\frac{nr}{2}=\begin{cases}
n(\lfloor\frac{n}{n-2m}\rfloor-1)/2,& \text{if $\lfloor\frac{n}{n-2m}\rfloor$ is odd;}\\
n(\lfloor\frac{n}{n-2m}\rfloor-2)/2,& \text{if   $\lfloor\frac{n}{n-2m}\rfloor$ is even.}
\end{cases}
\end{align*}




\section{Proofs of Theorems \ref{main-result1} and \ref{main-result2}}

For completing the proofs of    Theorems \ref{main-result1} and  \ref{main-result2}, we need the following result (see Exercise 3.3.18 in \cite{Lov}).

\begin{theorem}[Lov\'{a}sz and Plummer, \cite{Lov}]\label{GE-thm}
Let $G$ be a graph without perfect matchings. Then $G$ has a  subset $S$ of $V(G)$ such that
\begin{itemize}
\item [\rm{(i)}] every component of $G-S$ is factor-critical;
\item [\rm{(ii)}] there exists a maximum matching  in $G$ which matches each vertex of $S$ with vertices in different components of $G-S$.
    \end{itemize}
\end{theorem}

\noindent\textbf{Proof of Theorem \ref{main-result1}.}
If  $ n \leq 2m + 1$,  there is no  matching of size $m+1$ in any $n$-vertex graph, which implies  that  there does not exist $n$-vertex $(m+1)K_2$-saturated graphs, and thus  $rsat(n, (m+1)K_2)$ does not exist.

Next we consider $n\ge 2m+2$. If $n\le 3m$ and $m/(n-2m)$ is an integer, then we  can take a graph $G'$ consisting of $(n-2m)$ cliques $K_{\frac{2m}{n-2m}+1}$. Notice that $G'$ is a $\frac{2m}{n-2m}$-regular graph with order $n$ such that
 \begin{align*}
 \nu(G')=(n-2m)\frac{m}{n-2m}=m
 \end{align*}
 and
\begin{align*}
e(G')=\frac{n\frac{2m}{n-2m}}{2}=\frac{nm}{n-2m},
\end{align*}
and adding an edge between any two of $n-2m$ cliques to $G'$ yields a matching of size $m+1$. So the $n$-vertex regular $(m+1)K_2$-saturated graphs exist under this condition.

Assume that $G$ is an $r$-regular  $(m+1)K_2$-saturated graph of order $n$  such that   $e(G)=rsat(n,(m+1)K_2)$. 
 Since $n\ge 2m+2$, $G$ does not contain perfect matchings and $r\ge 2$. So by  Theorem \ref{GE-thm} (i), there is a subset $S$ of $V(G)$ such that every component of $G -S$ is factor-critical. We use  $D_1, D_2,\cdots,D_q$ to denote these  components of $G-S$. Let $s:=|S|$.  For $1\leq i\leq q$, write  $d_i:=|V(D_i)|$.  
By Theorem \ref{GE-thm} (i) and (ii), we can infer that
\begin{align}\label{eq1'}
\nu(G)=s+\sum_{i=1}^{q}\frac{d_i-1}{2}= m.
\end{align}

\begin{claim}\label{c1'}
$S = \emptyset$.
\end{claim}

By contradiction. Suppose that $S \neq \emptyset$. Then for any vertex $v \in S$,  we have $d_G(v) = n - 1$. Indeed, if  $d_G(v) < n - 1$,  there exists a vertex $u \in V(G) \setminus \{v\}$ such that $uv \notin E(G)$, and thus  by (\ref{eq1'}) we have   
\begin{align*}
\nu(G+\{uv\}) =s+\sum_{i=1}^{q}\frac{d_i-1}{2} = m=\nu(G),
\end{align*}
contradicting the assumption that $G$ is $(m+1)K_2$-saturated. 
Since $G$ is $r$-regular, it follows that $r = n - 1$. So $G$ is a complete graph with order $n\ge 2m+2$, which implies that there exists a  matching of size $m+1$ in $G$, again a contradiction. Hence, the claim holds.

\begin{claim}\label{c2'}
$d_1=d_2=\cdots=d_q\ge 3$.
\end{claim}

 Suppose, by contradiction, that there exist  two different components $D_i$ and $D_j$ with  $d_i < d_j$. Since $G$ is $r$-regular and $S=\emptyset$ by Claim \ref{c1'}, both $D_i$ and $D_j$ are $r$-regular. So we have $r \leq d_i - 1 < d_j - 1$, implying that $D_j$ is not a clique. Then there exists an edge $uv\notin E(D_j)$,  by (\ref{eq1'}),  we have 
\begin{align*}
\nu(G+\{uv\}) =s+\sum_{i=1}^{q}\frac{d_i-1}{2}= m =\nu(G),
\end{align*}
contradicting the assumption that $G$ is $(m+1)K_2$-saturated.  Therefore, all components $D_i$ have the same order. By Theorem \ref{GE-thm} (i), each $D_i$ is factor-critical, and so each has an odd order with at least $3$.

\begin{claim}\label{c3'}
For $1\le i\le q$, each $D_i$ is a clique.
\end{claim}

By contradiction. Suppose that $D_i$ is not a clique for some  $1\leq i\leq q$. Then adding an edge $uv$ within $D_i$ gives a graph $G+\{uv\}$, and by (\ref{eq1'}) we have 
\begin{align*}
\nu(G+\{uv\}) =s+\sum_{i=1}^{q}\frac{d_i-1}{2} = m=\nu(G),
\end{align*}
contradicting the $(m+1)K_2$-saturation of $G$. Therefore, the claim is ture.

By Claims \ref{c1'}, \ref{c2'}, and \ref{c3'}, we have that $G$ is the union of $q$ factor-critical cliques with the same order $d_1$. Thus
\begin{align}\label{eqn}
n=qd_1,
\end{align} 
and by (\ref{eq1'}), we have
\begin{align}\label{eqnu}
q\Big(\frac{d_1-1}{2}\Big)= m.
\end{align}
Combining  (\ref{eqn}) and (\ref{eqnu}), one can see 
\begin{align*}
(n-2m)(d_1-1)= q(d_1-1)= 2m,
\end{align*}
i.e.,
\begin{align*}
d_1-1= \frac{2m}{n-2m}.
\end{align*}
Therefore,  $2m/(n-2m)$ is an even integer (i.e., $m/(n-2m)$ is an integer), and 
\begin{align*}
e(G)=\frac{nr}{2}= \frac{n(d_1-1)}{2}= \frac{nm}{n-2m}.
\end{align*}

It remains to show that  $n\le 3m$.  Combining  (\ref{eqn}) and (\ref{eqnu}),  we obtain
\begin{align*}
q=n-2m.
\end{align*}
Thus, by Claim \ref{c2'},  we have
\begin{align*}
n= qd_1\ge 3(n-2m),
\end{align*}
i.e.,
$n\le 3m$.

The  proof of Theorem \ref{main-result1} is  complete. \qed

 \vspace{3ex} 
 
\noindent\textbf{Proof of Theorem \ref{main-result2}.}
We first show the lower bound for $rex(n,(m+1)K_2)$.
Let $r\in \{\lfloor n/n-2m\rfloor-1, \lfloor n/n-2m\rfloor-2\}$ be an even integer. Then we can see that  $n-(n-2m-1)(r+1)\geq r+1$ is odd through calculations. Consider the graph $G'$ consisting of $G_1\cup \bigcup_{i=2}^{n-2m}G_{i}$, where $G_i=K_{r+1}$ for each $i\in \{2,3,\ldots,n-2m\}$ and $G_1$ is an  $r$-regular factor-critical graph of order $n-(n-2m-1)(r+1)$. Notice that $G'$ is an $r$-regular graph of order $n$ with
 \begin{align*}
 \nu(G')=\frac{r}{2}(n-2m-1)+\frac{n-(n-2m-1)(r+1)-1}{2}=m,
 \end{align*}
and 
\begin{align*}
e(G')=\frac{nr}{2}= \begin{cases}
n(\lfloor\frac{n}{n-2m}\rfloor-1)/2,& \text{if  $\lfloor\frac{n}{n-2m}\rfloor$ is odd;}\\
n(\lfloor\frac{n}{n-2m}\rfloor-2)/2,& \text{if  $\lfloor\frac{n}{n-2m}\rfloor$ is even.}
\end{cases}
\end{align*}
Let  $G$ be a $d$-regular graph of order $n$  such that  $\nu(G)\leq m$ and $e(G)=rex(n,(m+1)K_2)$. By the definition of $rex(n,(m+1)K_2)$, 
\begin{align}\label{eq0}
e(G)\geq e(G')= \begin{cases}
n(\lfloor\frac{n}{n-2m}\rfloor-1)/2,& \text{if  $\lfloor\frac{n}{n-2m}\rfloor$ is odd;}\\
n(\lfloor\frac{n}{n-2m}\rfloor-2)/2,& \text{if  $\lfloor\frac{n}{n-2m}\rfloor$ is even.}
\end{cases}
\end{align}

In the following, we will prove $e(G)\leq e(G')$. Since $n\ge 2m+2$, $G$ does not contain perfect matchings and $d\ge 2$. So by  Theorem \ref{GE-thm} (i), there is a subset $S$ of $V(G)$ such that every component of $G -S$ is factor-critical. We use  $D_1, D_2,\cdots,D_q$ to denote these  components of $G-S$. Let $s:=|S|$.  For $1\leq i\leq q$, write  $d_i:=|V(D_i)|$.  
By Theorem \ref{GE-thm} (i) and (ii), we can infer that
\begin{align}\label{eq1}
m=s+\sum_{i=1}^{q}\frac{d_i-1}{2}
\end{align}
and
\begin{align}\label{eq2}
n=s+\sum_{i=1}^{q}d_i.
\end{align}
Combining (\ref{eq1}) and (\ref{eq2}), we have 
\begin{align}\label{eq3}
q=n-2m+s.
\end{align}

\begin{claim}\label{c1}
$d$ is even.
\end{claim}

Since the order of factor-critical components is odd, the claim holds clearly for the case of  $S=\emptyset$.  So we assume that $S\neq \emptyset$. We prove the claim by contradiction. Suppose that  $d$ is odd, then   $e_G(V(D_i),S)\geq 1$ for each $i\in\{1,2,\ldots, q\}$.  We may assume, without loss of generality,  that
\begin{align*}
d_1\geq \cdots \geq d_l\geq d+2>
d\geq d_{l+1}\geq \cdots \geq d_{q}\geq 1.
\end{align*}
For each $j\in \{l+1,\ldots, q\}$,  we have $d_{D_j}(v)\leq d_j-1$ for every vertex $v\in V(D_j)$.
Since $G$ is a $d$-regular graph,  for each $l+1\leq j \leq q$, we have 
\begin{align*}
e_G(V(D_j),S)\geq (d-(d_j-1))d_j\geq d,
\end{align*}
and
\begin{align*}
e_G\Big(S,\bigcup_{i=1}^{q}V(D_i)\Big)\leq ds.
\end{align*}
Thus, $q-l\leq s$, i.e., there are at most $s$ components of order less than $d+2$. Combining $q-l\le s$  and (\ref{eq3}), we have $l\geq n-2m$, i.e., there are at least $n-2m$ components of order (at least) $d+2$. One can see that
\begin{align*}
n\geq \Bigg|S\cup \bigcup_{i=1}^{l}V(D_i)\Bigg|\geq s+(n-2m)(d+2)>(n-2m)(d+2),
\end{align*} 
which implies that
\begin{align*}
d\leq
\begin{cases}
\lfloor\frac{n}{n-2m}\rfloor-2,& \text{if  $\lfloor\frac{n}{n-2m}\rfloor$ is odd;}\\
\lfloor\frac{n}{n-2m}\rfloor-3,& \text{if  $\lfloor\frac{n}{n-2m}\rfloor$ is even.}
\end{cases}
\end{align*}
Hence, 
\begin{align*}
e(G)=\frac{nd}{2}<\begin{cases}
\frac{n(\lfloor\frac{n}{n-2m}\rfloor-1)}{2},& \text{if  $\lfloor\frac{n}{n-2m}\rfloor$ is odd;}\\
\frac{n(\lfloor\frac{n}{n-2m}\rfloor-2)}{2},& \text{if  $\lfloor\frac{n}{n-2m}\rfloor$ is even,}
\end{cases}
\end{align*}
 contradicting with (\ref{eq0}).  This completes the proof of Claim \ref{c1}.

It follows from (\ref{eq1}) that  $0\leq s\leq m$. Two cases will be distinguished, and by Claim \ref{c1}, we   let $d$ be an  even integer.

\begin{case}
$1\leq s\leq m$.
\end{case}

 Let $D_1,D_2,\ldots,D_k$ be the components such that $E_{G}(V(D_i),S)\neq \emptyset$. Then each of the other  components $D_j\in \{D_{k+1},\ldots, D_{q}\}$ satisfies $e_{G}(V(D_j),S)= 0$ and $|V(D_j)|\geq d+1$ since $G$ is $d$-regular.
We may assume, without loss of generality,  that
\begin{align*}
d_1\geq \cdots \geq d_l\geq d+1
>
d-1\geq d_{l+1}\geq \cdots \geq d_{k}\geq 1.
\end{align*}
For each component  $D_{j}\in \{D_{l+1},\ldots, D_k\}$,  we have $d_{D_j}(v)\leq d_j-1$ for every vertex $v\in V(D_j)$.
Since $G$ is a $d$-regular graph,  for each $l+1\leq j \leq k$
\begin{align*}
e_G(V(D_j),S)\geq (d-(d_j-1))d_j\geq d,
\end{align*}
and
\begin{align*}
e_G\Big(S,\bigcup_{i=1}^{k}V(D_i)\Big)\leq ds.
\end{align*}
Thus, $k-l\leq s$, i.e., there are at most $s$ components of order less than $d+1$. By (\ref{eq3}), $l+(q-k)\geq n-2m$, i.e., there are at least $n-2m$ components of order (at least) $d+1$. One can see that
\begin{align*}
n&\geq \Bigg|S\cup \bigcup_{i=1}^{l}V(D_i)\Bigg|\\
&\geq s+(n-2m)(d+1)\\
&>(n-2m)(d+1)\\
&\ge 3n-6m~~(since~d\ge 2),
\end{align*} 
which implies that $n<3m$ and 
\begin{align*}
d\leq
\begin{cases}
\lfloor\frac{n}{n-2m}\rfloor-1,& \text{if  $\lfloor\frac{n}{n-2m}\rfloor$ is odd;}\\
\lfloor\frac{n}{n-2m}\rfloor-2,& \text{if  $\lfloor\frac{n}{n-2m}\rfloor$ is even.}
\end{cases}
\end{align*}
Hence, 
\begin{align*}
e(G)=\frac{nd}{2}\leq\begin{cases}
\frac{n(\lfloor\frac{n}{n-2m}\rfloor-1)}{2},& \text{if  $\lfloor\frac{n}{n-2m}\rfloor$ is odd;}\\
\frac{n(\lfloor\frac{n}{n-2m}\rfloor-2)}{2},& \text{if  $\lfloor\frac{n}{n-2m}\rfloor$ is even.}
\end{cases}
\end{align*}

\begin{case}
$s=0$.
\end{case}

Note that each component $D_i\in \{D_1,D_2,\ldots,D_q\}$ has at least $d+1$vertices since $G$ is $d$-regular. Then by (\ref{eq3}), we have 
\begin{align*}
n\geq q(d+1)=(n-2m)(d+1)\ge 3n-6m,
\end{align*}
which implies that $n\le 3m$ and 
\begin{align*}
d\leq
\begin{cases}
\lfloor\frac{n}{n-2m}\rfloor-1,& \text{if  $\lfloor\frac{n}{n-2m}\rfloor$ is odd;}\\
\lfloor\frac{n}{n-2m}\rfloor-2,& \text{if  $\lfloor\frac{n}{n-2m}\rfloor$ is even.}
\end{cases}
\end{align*}
Hence, 
\begin{align*}
e(G)=\frac{nd}{2}\leq\begin{cases}
\frac{n(\lfloor\frac{n}{n-2m}\rfloor-1)}{2},& \text{if  $\lfloor\frac{n}{n-2m}\rfloor$ is odd;}\\
\frac{n(\lfloor\frac{n}{n-2m}\rfloor-2)}{2},& \text{if  $\lfloor\frac{n}{n-2m}\rfloor$ is even.}
\end{cases}
\end{align*}

The proof of Theorem \ref{main-result2} is complete. \qed

\section*{Declaration of competing interest}

We declare that we have no financial and personal relationships with other people or organizations that can inappropriately influence our work, there is no professional or other personal interest of any nature or kind in any product, service and/or company that could be construed as influencing the position presented in, or the review of, the manuscript entitled.

\section*{Acknowledgement}
The research of the first author was supported by JST SPRING (No. JPMJSP2178). The research of the second author was supported by National Key Research and Development Program of China (No.~2023YFA1010203), and National Natural Science Foundation of China (Nos.~12401464 and 12471334).  The research of the third author was supported by National Natural Science Foundation of China (Nos.~12131013 and 12471334) and Shaanxi Fundamental Science Research Project for Mathematics and Physics (No. 22JSZ009).

\section*{Data availability}

No data was used for the research described in the article.


\end{document}